\documentclass[11pt,a4paper]{article}

\usepackage{epsf,epsfig,amsfonts,amsgen,amsmath,amstext,amsbsy,amsopn,amsthm
}
\usepackage{enumerate}
\usepackage{amsmath}
\usepackage{amsfonts,amsthm,amssymb}
\usepackage{amsfonts}
\usepackage{graphics}
\usepackage{latexsym,bm}
\usepackage{amsfonts,amsthm,amssymb,bbding}
\usepackage{indentfirst}
\usepackage{graphicx}
\usepackage{color}
\usepackage{moresize}
\usepackage{array}
\usepackage{multirow}
\usepackage{graphicx}
\usepackage{threeparttable}
\usepackage{makecell}
\renewcommand\proofname{\bf Proof}
\newtheorem{theorem}{Theorem}%[section]

\newtheorem{lemma}{Lemma}[section]

\newtheorem{false statement}{False statement}

\theoremstyle{definition}

\newtheorem{claim}{Claim}

\newtheorem{remark}[claim]{Remark}

\begin{document}

\title{\bf\Large Toughness and spectral radius in graphs \thanks{This work is supported by the Scientific Research Plan of Universities in Xinjiang, China (No. XJEDU2022P009), the National Natural Science Foundation of China (Grant No. 12271162 and 12301454), the Natural Science Foundation of Shanghai (No. 22ZR1416300) and The Program for Professor of Special Appointment (Eastern Scholar) at Shanghai Institutions of Higher Learning (No. TP2022031).}}
\author{\Large Yuanyuan Chen$^{a}$, Dandan Fan$^{b,c}$\thanks{Corresponding author;
Email addresses:
chenyy\underline{\ }de@sina.com (Y. Chen),ddfan0526@163.com (D. Fan), huiqiulin@126.com (H. Lin).}, Huiqiu Lin$^{b}$,\\
\small\it $^{a}$ College of Mathematics and System Science, Xinjiang University,\\
\small\it  Urumqi, Xinjiang 830017, China\\[1mm]
\small\it $^b$ School of Mathematics, East China University of Science and Technology, \\
\small\it   Shanghai 200237, China\\[1mm]
\small\it $^c$ College of Mathematics and Physics, Xinjiang Agricultural University,\\
\small\it Urumqi, Xinjiang 830052, China\\[1mm]
}
\date{}

\maketitle

\begin{center}
\begin{minipage}{140mm}
\begin{center}{\bf Abstract}\end{center}
{

The Brouwer's toughness conjecture states that every $d$-regular connected graph always has $t(G)>\frac{d}{\lambda}-1$ where $\lambda$ is the second largest absolute eigenvalue of the adjacency matrix.
In 1988,  Enomoto introduced a variation of toughness $\tau(G)$ of a graph $G$. By incorporating the variation of toughness and spectral conditions,
we provide spectral conditions for a graph to be $\tau$-tough ($\tau\geq 2$ is an integer)
and to be $\tau$-tough ($\frac{1}{\tau}$ is a positive integer) with minimum degree $\delta$, respectively.
Additionally, we also investigate a analogous problem concerning balanced bipartite graphs.

}
\end{minipage}

\end{center}

\renewcommand{\thefootnote}{}
%\footnote{E-mail address: jimingguo@hotmail.com (J.-M. Guo)}
\vskip 0.3in

\section{Introduction}
Throughout this paper, we consider only finite, undirected and simple connected graphs.
The \textit{toughness $t(G)$} of a non-complete graph $G$ is defined as $t(G)=\min\Big\{\frac{|S|}{c(G-S)}\Big\}$ in which the minimum is taken over all proper sets $S\subset G$ where $c(G-S)$ denotes the number of components of $G-S$. The concept of toughness was introduced by Chv\'atal \cite{V. C} in 1973 to capture combinatorial properties related to the cycle structure of a graph. The toughness is related to many other important properties of a graph, such as the existence of various factors \cite{H. Enomoto, H. Enomoto1}, cycles \cite{M. Fiedler} and spanning trees \cite{D. Bauer1}. For more extensive work on toughness, one can see the survey \cite{D. Bauer2}.

%and the existence of various factors \cite{H. Enomoto, H. Enomoto1}, cycles or spanning trees \cite{D. Bauer1}.
%Toughness is a crucial parameter considered in network security which characterizes
%the vulnerability of the network from the perspective of graph topology and it is a hard parameter to determine exactly.
%The toughness of a graph is related to many other important properties of a graph such as Hamiltonicity \cite{M. Fiedler, D. Fan2}
%and the existence of various factors \cite{H. Enomoto, H. Enomoto1}, cycles or spanning trees \cite{D. Bauer1}.
%Toughness and special graph classes, for example triangle-free graphs have a number of interesting properties
%\cite{N. Alon, D. Bauer3, T. Bohme} etc..
%Enomoto et al. \cite{H. Enomoto, H. Enomoto1} proved Chv\'atal's \cite{V. C} conjecture that $k$-tough graphs have
%$k$-factors if they satisfy trivial necessary conditions.

The study of the relationship between toughness on eigenvalues was initiated by Alon \cite{Alon} who showed that for any connected $d$-regular graph $G$, $t(G)>\frac{1}{3}(\frac{d^2}{d\lambda+\lambda^2}-1)$ where $\lambda$ is the second largest absolute eigenvalue of the adjacency matrix. Meanwhile, Brouwer \cite{Brouwer1} independently proved $t(G)>\frac{d}{\lambda}-2$, and he \cite{Brouwer2} further conjectured that $t(G)>\frac{d}{\lambda}-1$. Recently, Gu \cite{Gu2} confirmed this conjecture completely. Very recently, Fan, Lin and Lu extended the results in terms of its spectral radius, and provided spectral conditions for a graph with minimum degree to be 1-tough and to be $t$-tough, respectively. For more details, we refer the reader to \cite{CG,CW,GH}.

In order to better investigate the existence of factors in a graph, Enomoto \cite{H. Enomoto2} introduced a variation of toughness in 1998. A non-complete graph $G$ is \textit{$\tau$-tough} if $|S|\geq \tau(c(G-S)-1)$ for every proper subset $S\subseteq V(G)$ with $c(G-S)>1$. The \textit{variation of toughness $\tau(G)$} of $G$ is the maximum $\tau$ for which $G$ is $\tau$-tough. In this paper, we consider the case where $|S|$ and $c(G-S)-1$ are divisible by each other.
Inspired by the work of Fan, Lin and Lu \cite{D. Fan2}, it is interesting to find spectral radius conditions for a graph to be $\tau$-tough, where $\tau$ or $\frac{1}{\tau}$ is a positive integer.
%Let $\mathcal{G}_{n,\tau}$ be the set of connected graphs of order $n$ which are not $\tau$-tough, where $\tau$ or
%$\frac{1}{\tau}$ is a positive real number.

The largest eigenvalue of $A(G)$, denoted by $\rho(G)$, is called the \emph{spectral radius} of $G$. Given two graphs $G_1$ and $G_2$, the \textit{disjoint union} $G_1\cup G_2$ is the graph with vertex set $V(G_1)\cup V(G_2)$ and edge set $E(G_1)\cup E(G_2)$, and the \textit{join}  $G_{1} \nabla G_{2}$ is the graph obtained from  $G_{1}\cup G_{2}$ by adding all edges between $G_{1}$ and $G_{2}$. Let $\delta(G)$ be the minimum degree of $G$.

%The spectrum of a graph is related to many important combinatorial parameters.
%The relationship between the toughness of a regular graph and the eigenvalues has
%been considered by many researchers, see \cite{N. Alon, A.E. B, S.M. C1, S.M. C2}.
%For positive integers $a\le b$, an even (or odd) $[a,b]$-factor of a graph $G$ is a
%spanning subgraph $H$ such that for each vertex $v\in V(G)$, $d_H(v)$ is even (or odd)
%and $a\le d_H(v)\le b$.
%Fan, Lin and Lu \cite{D. Fan1} provided spectral conditions for the existence of an odd  $[1,b]$-factor
%in a connected graph with minimum degree $\delta$.
%In this paper, we focus on the variation of toughness $\tau(G)$.
%We first provide a spectral condition for a graph to be $\tau$-tough ($\tau$ is a positive integer) and to be $\frac{1}{b}$-tough
%($b$ is a positive integer) with minimum degree $\delta$, respectively.
%From the well know sufficient and necessary condition for the existence of an odd $[1,b]$-factor was established by Amahashi \cite{A. Amahashi},
%we get that Fan, Lin and Lu's [\cite{D. Fan1}, Theorem 1.2] conclusion is a corollary of our results, and our results can improve their conclusion.

\begin{theorem}\label{ta}
Suppose that $G$ is a connected graph of order $n$. Thus, the following statements hold.
\begin{enumerate}[(i)]
\item Let $\tau\ge2$ be an integer and $n\ge 2\tau^2+3\tau$. If $\rho(G)\ge \rho(K_{\tau-1}\nabla(K_{n-\tau}\cup K_1))$, then $G$ is a $\tau$-tough graph, unless $G\cong K_{\tau-1}\nabla(K_{n-\tau}\cup K_1)$.
\item Let $\frac{1}{\tau}\geq 1$ and $\delta$ be two positive integers, and let $n\ge\max\{5\delta+4, \frac{\delta^3}{\tau}+\delta\}$. If $\delta(G)=\delta $ and $\rho(G)\ge \rho\Big(K_{\delta}\nabla \Big(K_{n-\frac{(\tau+1)\delta}{\tau}-1}\cup \Big(\frac{\delta}{\tau}+1\Big)K_1\Big)\Big),$ then
$G$ is a $\tau$-tough graph, unless $G\cong K_{\delta}\nabla$ $\Big(K_{n-\frac{(\tau+1)\delta}{\tau}-1}\cup (\frac{\delta}{\tau}+1)K_1\Big)$.
\end{enumerate}
\end{theorem}

One can verify that $t(G)$ at most 1 for any bipartite graph $G=(X,Y)$ because removing $X$ from $G$ (assuming $|X|\ge|Y|$) results in an independent set $Y$. In 2006, Bian \cite{Q. Bian} defined the toughness of a bipartite graph with the aim of providing a better description of its toughness. The \emph{bipartite toughness} $t^B(G)$ of a non-complete bipartite graph $G=(X,Y)$ is defined as $t^B(G)=\min\Big\{\frac{|S|}{c(G-S)}\Big\}$, in which the minimum is taken over all proper subsets $S\subset X$ (or $S\subset Y$) such that $G-S$ is disconnected and $c(G-S)>1$. A bipartite graph $G = (X, Y)$ is called \emph{balanced} if $|X| = |Y|$. In the past few decades, many researchers utilized the toughness $t^B(G)$ to investigate the existence of factors in balanced bipartite graphs \cite{P. Katerinis,D. Kratsch,Liu, Wang}. Analogous to the toughness of bipartite graph, we introduce a variation of bipartite toughness. A balanced non-complete bipartite graph $G=(X,Y)$ is \textit{$\tau^B$-tough} if $|S|\geq \tau^B(c(G-S)-1)$ for every proper subset $S\subset X$ (or $S\subset Y$) with $c(G-S)> 1$. The \textit{variation of bipartite toughness $\tau^B(G)$} of a  balanced non-complete bipartite $G$ is the maximum $\tau^B$ for which $G$ is $\tau^B$-tough. In this paper, we also provide spectral conditions for a balanced bipartite graph to be $\tau^B$-tough where $\tau^B$ or $\frac{1}{\tau^B}$ is a positive integer.

Given two bipartite graphs $G_1=(X_{1}, Y_{1})$
and $G_2=(X_{2}, Y_{2})$, let $G_{1}\nabla_{1} G_2$ denote the graph obtained from $G_1 \cup G_2$ by adding all possible
edges between $X_1$ and $Y_2$, $X_{2}$ and $Y_{1}$. For every two subsets $A$ and $B$ of $V(G)$,
let $e(A,B)$ denote the number of edges with one end in $A$ and the other one in $B$.
Let $O_{a,b}$ be a bipartite graph with $|A|=a$, $|B|=b$ and $e(A,B)=0$.
Suppose that $r\geq 2$ be an integer. For  $2r\mid n$, let
$B_n^r=K_{\frac{n}{2}-1,\frac{n}{2}-\frac{n}{2r}}\nabla_1 O_{1,\frac{n}{2r}}$. For $2r\nmid n$, let
 $B_n^r\in\{K_{r\lfloor\frac{n}{2r}\rfloor-1,\frac{n}{2}-\lfloor\frac{n}{2r}\rfloor}\nabla_1 O_{\frac{n}{2}-r\lfloor\frac{n}{2r}\rfloor+1,\lfloor\frac{n}{2r}\rfloor},K_{r-1,\frac{n}{2}-1}\nabla_1 O_{\frac{n}{2}-r+1,1}\}$ with $\rho(B_n^r)=\max\{\rho(K_{r\lfloor\frac{n}{2r}\rfloor-1,\frac{n}{2}-\lfloor\frac{n}{2r}\rfloor}\nabla_1 O_{\frac{n}{2}-r\lfloor\frac{n}{2r}\rfloor+1,\lfloor\frac{n}{2r}\rfloor}),\rho(K_{r-1,\frac{n}{2}-1}\nabla_1 O_{\frac{n}{2}-r+1,1})\}$.
%$$
%B_n^r\cong \left\{
%\begin{array}{ll}
%K_{r\lfloor\frac{n}{2r}\rfloor-1,\frac{n}{2}-\lfloor\frac{n}{2r}\rfloor}\nabla_1 O_{\frac{n}{2}-r\lfloor\frac{n}{2r}\rfloor+1,\lfloor\frac{n}{2r}\rfloor} & \mbox{if $\rho_1\leq \rho_2$,}\vspace{1ex}\\
%K_{r-1,\frac{n}{2}-1}\nabla_1 O_{\frac{n}{2}-r+1,1} & \mbox{if $\rho_1< \rho_2$ }
%\end{array}
%\right.
%$$

\begin{theorem}\label{td}
Suppose that $\tau^B=r$ and $G$ is a connected balanced bipartite graph of order $n$. Thus, the following statements hold.
\begin{enumerate}[(i)]
\item Let $r\geq 2$ be an integer and $n\geq 2r^2+6r$. If
$\rho(G)\ge\rho(B_n^r),$
then $G$ is a $r$-tough graph, unless $G\cong B_n^r.$
\item Let $\frac{1}{r}$ be a positive integer and $n\geq \frac{4}{r}+6$. If
$\rho(G)\ge \rho(K_{1,\frac{n}{2}-\frac{1}{r}-1}\nabla_1 O_{\frac{n}{2}-1,\frac{1}{r}+1}),$
then $G$ is a $r$-tough graph, unless $G\cong K_{1,\frac{n}{2}-\frac{1}{r}-1}\nabla_1 O_{\frac{n}{2}-1,\frac{1}{r}+1}$.
\end{enumerate}
\end{theorem}

\begin{remark}
Let $G=K_{r\lfloor\frac{n}{2r}\rfloor-1,\frac{n}{2}-\lfloor\frac{n}{2r}\rfloor}\nabla_1 O_{\frac{n}{2}-r\lfloor\frac{n}{2r}\rfloor+1,\lfloor\frac{n}{2r}\rfloor}$ and $G'=K_{r-1,\frac{n}{2}-1}\nabla_1 O_{\frac{n}{2}-r+1,1}$. The comparison of the spectral radius of these two graphs is not always straightforward. For example, by Matlab programming, we have $18.499\approx \rho(G')>\rho(G)\approx 18.472$ if $r=3$ and $ n=38$, and $134.50\approx \rho(G')>\rho(G)\approx 134.46$ if $r=10$ and $ n=270$. However, for  $r=10$ and $n=402$, we get $200.50\approx \rho(G')<\rho(G)\approx 200.81$.
\end{remark}

\section{Preliminaries}

Let $M$ be a real $n\times n$ matrix, and let $X=\{1,2,\ldots,n\}$. Given a partition $\Pi:X=X_{1}\cup X_{2}\cup \cdots \cup X_{k}$,  the matrix $M$ can be correspondingly partitioned as
$$
M=\left(\begin{array}{ccccccc}
M_{1,1}&M_{1,2}&\cdots &M_{1,k}\\
M_{2,1}&M_{2,2}&\cdots &M_{2,k}\\
\vdots& \vdots& \ddots& \vdots\\
M_{k,1}&M_{k,2}&\cdots &M_{k,k}\\
\end{array}\right).
$$
The \textit{quotient matrix} of $M$ with respect to $\Pi$ is defined as the $k\times k$ matrix $B_\Pi=(b_{i,j})_{i,j=1}^k$ where $b_{i,j}$ is the  average value of all row sums of $M_{i,j}$.
The partition $\Pi$ is called \textit{equitable} if each block $M_{i,j}$ of $M$ has constant row sum $b_{i,j}$.
Also, we say that the quotient matrix $B_\Pi$ is \textit{equitable} if $\Pi$ is an equitable partition of $M$.

\begin{lemma}(Brouwer and Haemers \cite[p. 30]{BH}; Godsil and Royle \cite[pp. 196--198]{C.Godsil})\label{Q}
Let $M$ be a real symmetric matrix, and let $\lambda_{1}(M)$ be the largest eigenvalue of $M$. If $B_\Pi$ is an equitable quotient matrix of $M$, then the eigenvalues of  $B_\Pi$ are also eigenvalues of $M$. Furthermore, if $M$ is nonnegative and irreducible, then $\lambda_{1}(M) = \lambda_{1}(B_\Pi).$
\end{lemma}

By the well-known Perron-Frobenius theorem (cf. \cite[Section 8.8]{C.Godsil}), we can easily deduce the following result.

\begin{lemma}\label{edge}
If $H$ is a spanning subgraph of a connected graph $G$, then
$$\rho(H)\leq\rho(G),$$
with equality if and only if $H\cong G$.
\end{lemma}

\section{Proof of Theorem \ref{ta}}
%Let $A$ be a real symmetric matrix of order $n$, and let $X=\{1,2,\ldots,n\}$.
%For any partition $\Pi$ of $X:$ $X_1\cup X_2\cup \cdots \cup X_m$, the matrix $A$ can be partitioned correspondingly.
%\begin{align*}
%A=\left[
%\begin{matrix}
%  A_{1,1} &\cdots & A_{1,m} \\
%  \vdots& \ddots  & \vdots \\
%  A_{m,1} &\cdots & A_{m,m}
%\end{matrix}
%  \right].
%\end{align*}
%The \emph{characteristic matrix} of $\Pi$ is the $n\times m$ matrix $\chi_\Pi$ whose columns are
%the characteristic vectors of $X_1,\ldots, X_m$, and \emph{quotient matrix} of $A$ with respect to $\Pi$
%is the matrix $B_\Pi=(b_{i,j})_{m\times m}$ with $b_{i,j}=\frac{1}{|X_i|}e^T_{|X_i|}A_{i,j}e_{|X_j|}$,
%where $e_{|X_i|}$ and $e_{|X_j|}$ are the all ones $|X_i|$- and $|X_j|$-vectors, respectively.
%In particular, the partition $\Pi$ is called \emph{equitable} if each block $A_{i,j}$ has constant row sum.
%Let $B_{\Pi}$ be a quotient matrix of $A$ with respect to some partition $\Pi$. From {\rm\cite{A.E., W.H.}}
%we know the following facts. The eigenvalues of
%$B_{\Pi}$ interlace those of $A$. Furthermore, if the partition $\Pi$ is equatable, then all eigenvalues of
%$B_{\Pi}$ are also eigenvalues of $A$.

%Let $G$ be a graph with adjacency matrix $A(G)$. The largest eigenvalue of $A(G)$, denoted by $\rho(G)$, is called the spectral
%radius of $G$. Denote by $\nabla$ and $\cup$ the join and union products, respectively.
\begin{lemma}[See \label{Dan Fan}\cite{D. Fan3}]
Let $n=\sum_{i=1}^tn_i+s$. If $n_1\ge n_2\ge \cdots\ge n_t\ge p$ and $n_1<n-s-p(t-1)$, then
$$\rho(K_s\nabla(K_{n_1}\cup K_{n_2}\cup\cdots\cup K_{n_t}))<\rho(K_s\nabla(K_{n-s-p(t-1)}\cup (t-1)K_p)).$$
\end{lemma}

\begin{lemma}[See \label{Hong}\cite{Y. Hong2}]
Let $G$ be a graph on $n$ vertices and $m$ edges. Then
$$\rho(G)\le\sqrt{2m-n+1},$$
with equality if and only if $G$ is a star or a complete graph.
\end{lemma}

\begin{lemma}[See \label{degree}\cite{Y. Hong1, V. Nikiforov1}]
Let $G$ be a graph on $n$ vertices and $m$ edges with minimum degree $\delta\ge1$. Then
$$\rho(G)\le\frac{\delta-1}{2}+\sqrt{2m-n\delta+\frac{(\delta+1)^2}{4}},$$
with equality if and only if $G$ is either a $\delta$-regular graph or a bidegreed graph in which each vertex is
of degree either $\delta$ or $n-1$.
\end{lemma}

%\begin{lemma}[See \cite{V. Nikiforov1}] \label{V. Nikiforov1}
%For nonnegative integers $p$ and $q$ with $2q \leq p(p-1)$ and $0 \leq x \leq p-1$, the function $f(x)=(x-%1) / 2+\sqrt{2 q-p x+(1+x)^{2} / 4}$ is decreasing with respect to $x$.
%\end{lemma}

 Let $e(G)$ denote the number of edges in $G$.
\renewcommand\proofname{\bf Proof of Theorem \ref{ta}}
\begin{proof}

Suppose that $G$ is not a $\tau$-tough graph with $\tau\geq 2$ or $\frac{1}{\tau}$ is a positive integer. Thus, there exists some nonempty subset $S$ of $V(G)$
such that $|S|<\tau(c(G-S)-1)$. Let $|S|=s$ and $c(G-S)=q$.

(i) Since $\tau\geq 2$ is an integer and $s\leq\tau(q-1)-1$, $G$ is a spanning
subgraph of $G'=K_{\tau(q-1)-1}\nabla(K_{n_1}\cup K_{n_2}\cup\cdots\cup K_{n_q})$
for some integers $n_1\ge n_2\ge\cdots\ge n_q$ with $\sum_{i=1}^qn_i= n-\tau(q-1)+1$.
Combining this with Lemma \ref{edge}, we have
\begin{equation}\label{equ::1}
\begin{aligned}
\rho(G)\le\rho(G'),
\end{aligned}
\end{equation}
with equality if and only if $G\cong G'$.
Let $G''=K_{\tau(q-1)-1}\nabla(K_{n-(\tau+1)(q-1)+1}\cup(q-1)K_1)$.
Then, by Lemma \ref{Dan Fan}, we have
\begin{equation}\label{equ::2}
\begin{aligned}
\rho(G')\le\rho(G''),
\end{aligned}
\end{equation}
with equality if and only if $G'\cong G''$. For the case of $q=2$, we have $G''\cong K_{\tau-1}\nabla(K_{n-\tau}\cup K_1)$. Combining this with (\ref{equ::1}) and (\ref{equ::2}), we deduce that
$$\rho(G)\leq \rho(K_{\tau-1}\nabla(K_{n-\tau}\cup K_1)),$$
with equality if and only if
$G\cong K_{\tau-1}\nabla(K_{n-\tau}\cup K_1)$.
For the case of $q\ge3$, according to Lemma \ref{Hong}, we have
\begin{equation}\label{equ::3}
\begin{aligned}
\rho(G'')&\le\sqrt{2e(G'')-n+1}\\
&=\sqrt{(n-q+1)(n-q)+2(\tau(q-1)-1)(q-1)-n+1}\\
         &=\sqrt{(1+2\tau)q^2-(2n+3+4\tau)q+n^2+2\tau+3}.
\end{aligned}
\end{equation}
Let $f(q)=(1+2\tau)q^2-(2n+3+4\tau)q+n^2+2\tau+3.$
Observe that $n\ge(\tau+1)(q-1)$. Thus $3\le q\le \frac{n}{\tau+1}+1$.
By a simple calculation, we have
 $$f(3)-(n-2)^2=-2n+8\tau-1<0$$
due to $n\geq 2\tau^2+3\tau$ and $\tau\geq 2$, and
\begin{equation*}
\begin{aligned}
f\Big(\frac{n}{\tau+1}+1\Big)-(n-2)^2&=-\frac{n^2-(2\tau^2+3\tau+1) n+3(\tau+1)^2}{(\tau+1)^2}\\
         &\le-\frac{\tau^2+3\tau+3}{(\tau+1)^2}~~(\mbox{since $n\geq 2\tau^2+3\tau$})\\
         &<0 ~~(\mbox{since $\tau\geq 2$}),
\end{aligned}
\end{equation*}
from which we get
$$f(q)\leq \max\Big\{f(3),f\Big(\frac{n}{\tau+1}+1\Big)\Big\}< (n-2)^2$$
for $3\le q\le \frac{n}{\tau+1}+1$. According to  (\ref{equ::3}), we can obtain that
\begin{equation}\label{equ::4}
\begin{aligned}
\rho(G'')\leq \sqrt{f(q)}<n-2.
\end{aligned}
\end{equation}
Since $K_{\tau-1}\nabla(K_{n-\tau}\cup K_1)$ contains $K_{n-1}$ as a proper subgraph,
we have $$\rho(K_{\tau-1}\nabla(K_{n-\tau}\cup K_1))>\rho(K_{n-1})=n-2$$ by Lemma \ref{edge}.
Combining this with (\ref{equ::1}), (\ref{equ::2}) and (\ref{equ::4}), we have
$$\rho(G)\le \rho(G')\le\rho(G'')<\rho(K_{\tau-1}\nabla(K_{n-\tau}\cup K_1)).$$

When  $\tau\geq 2$ is an integer, by the discussions as above, we have
$$\rho(G)\le \rho(K_{\tau-1}\nabla(K_{n-\tau}\cup K_1)),$$
with equality if and only if $G\cong K_{\tau-1}\nabla(K_{n-\tau}\cup K_1)$. As $K_{\tau-1}\nabla(K_{n-\tau}\cup K_1)$ is not a $\tau$-tough graph, we complete the proof of (i).

(ii) %Suppose that $G$ is not a $\tau$-tough graph. Thus, there exists some nonempty subset $S$ of $V(G)$
%such that $|S|<\tau(c(G-S)-1)$. Take $|S|=s$ and $c(G-S)=q$.
Notice that $\frac{1}{\tau}\geq 1$ is an integer. Thus, we may assume that $b=\frac{1}{\tau}$ ($b\geq 1$), and hence $q\ge bs+2$. One can verify that $G$ is a spanning subgraph of $G^1_s=K_s\nabla(K_{c_1}\cup K_{c_2}\cup\cdots\cup K_{c_l})$ for some
$c_1\ge c_2\ge\cdots\ge c_l$ with $l=bs+2$ and $\sum_{i=1}^lc_i=n-s$. Thus, by Lemma \ref{edge}, we have
\begin{equation}\label{equ::5}
\begin{aligned}
\rho(G)\le\rho(G^1_s),
\end{aligned}
\end{equation}
with equality if and only if $G\cong G^1_s$.
Then we shall divide the proof into the following three cases.

\noindent{\bf{Case 1.}} $s\ge \delta+1$.

Let $G^2_s=K_s\nabla(K_{n-(b+1)s-1}\cup (bs+1)K_1)$. By Lemma \ref{Dan Fan}, we have
\begin{equation}\label{equ::6}
\begin{aligned}
\rho(G_s^1)\le\rho(G^2_s),
\end{aligned}
\end{equation}
with equality if and only if $G\cong G^2_s$. Notice that $A(G^2_s)$ has the equitable quotient matrix
\begin{align*}
A_\Pi^s=\left[
\begin{matrix}
  s-1 & bs+1 & n-(b+1)s-1 \\
  s & 0 & 0  \\
  s & 0 & n-(b+1)s-2 \\
\end{matrix}
  \right].
\end{align*}
By a simple calculation, the characteristic polynomial of $A_\Pi^s$ is
$$\phi(A_\Pi^s,x)=x^3-(n-bs-3)x^2-(n+s+bs^2-bs-2)x-s(bs+1)(bs-n+s+2).$$
Note that $A(K_\delta\nabla(K_{n-(b+1)\delta-1}\cup (b\delta+1)K_1))$ has the equitable quotient matrix $A^\delta_\Pi$,
which is obtained by replacing $s$ with $\delta$ in $A_\Pi^s$. Then
\begin{equation}\label{equ::7}
\phi(A_\Pi^s,x)-\phi(A_\Pi^\delta,x)=(s-\delta)\varphi(x)
\end{equation}
where $$\varphi(x)=bx^2-((\delta+s-1)b+1)x-(\delta^2+\delta s+s^2)b^2-(\delta^2+(s+3-n)\delta+s^2-(n-3)s)b+n-\delta-s-2.$$
 The symmetry axis of parabola $\varphi(x)$ is
$$x=\frac{\delta+s-1}{2}+\frac{1}{2b}<n-b\delta-2$$
due to $n\ge(b+1)s+2$, $s\geq \delta+1$ and $b\geq 1$. This implies that $\varphi(x)$ is increasing with respect to $x\ge n-b\delta-2$. Since $n\ge(b+1)s+2$, we have $s\le\frac{n-2}{b+1}$. Hence
\begin{equation}\label{equ::8}
\begin{aligned}
\varphi(x)&\ge\varphi(n-b\delta-2)\\
&=\!-\!(b^2\!+\!b)s^2\!-\!((\delta\! +\! 1)b\!+\!1)s\!+\!3b^2\delta\!-\!(\delta^2\!-\!2)b\!+\!bn^2\!-\!(2b^2\delta\!+\!3b)n\!+\!b^3\delta^2\!-\!\delta\\
&\ge\frac{1}{b+1}g(n) ~~(\mbox{since $s\le\frac{n-2}{b+1}$}),
\end{aligned}
\end{equation}
where $$g(n)=b^2n^2-(2b^3\delta+2b^2\delta+3b^2+b\delta+1)n+b\delta(b^3\delta+b^2\delta+3b^2-b\delta+3b-\delta+1)+2b^2-\delta+2.$$
If $b=1$, since $n\geq 5\delta+4$ and $\delta\geq 1$, we have $$g(n)=n(n-5\delta-4)+6\delta+4>0.$$
In what follows, we will consider the case of $b\geq 2$. If $\delta=1$, then
$$g(n)=b^2n^2-(2b^3+5b^2+b+1)n+b(b^3+4b^2+2b)+2b^2+1.$$
Observe that $s\geq\delta+1=2$ and $n\geq (b+1)s+2\geq 2b+4$. The symmetry axis of parabola $g(n)$ is
$$\frac{2b^3+5b^2+b+1}{2b^2}=b+\frac{5}{2}+\frac{1}{2b}+\frac{1}{2b^2}<2b+4$$
due to $b\geq 2$. Then $g(n)$ is increasing with respect to $n\geq 2b+4$, and hence $$g(n)\ge g(2b+4)=(b^2-3)(b+1)^2>0$$
due to $b\geq 2$.
If $\delta\geq 2$, since the symmetry axis of parabola $g(n)$ is
$$b\delta+\delta+\frac{\delta}{2b}+\frac{1}{2b^2}+\frac{3}{2}<b\delta^3+\delta,$$
we conclude that $g(n)$ is increasing with respect to $n\geq b\delta^3+\delta$.
 For $n\geq b\delta^3+\delta$, we get
\begin{equation*}
\begin{aligned}
  &g(n)\ge g(b\delta^3+\delta)\\
&=\delta^3b^2(\delta^3b^2-2\delta b^2-\delta-3b)+((b^2-b-2)\delta^2+2+3\delta b)b^2-(\delta^3+2\delta^2)b+\delta(b-2)+2\\  &\ge \delta^3b^2(2\delta b^2-\delta-3b)-(\delta^3+2\delta^2)b \\
  &=\delta^3b^2(2(\delta-2)(b-2)+2\delta b(b-1)+3\delta+b-8)-(\delta^3+2\delta^2)b\\
  &\ge 2\delta^4b^3-(\delta^3+2\delta^2)b\\
&=\delta^2b(2\delta^2b^2-\delta-2)\\
  &>0,
\end{aligned}
\end{equation*}
where all the inequalities follows from the fact that $\delta\ge2$ and $b\ge2$. According to the above discussions, we have $g(n)>0$ for $b\geq 1$ and $\delta\geq 1$.
Combining this with (\ref{equ::7}) and (\ref{equ::8}), we have $\phi(A_\Pi^s,x)>\phi(A_\Pi^\delta,x)$ for $x>n-b\delta-2$.
Since $K_{n-b\delta-1}$ is a proper subgraph of $ K_{\delta}\nabla (K_{n-(b+1)\delta-1}\cup (b\delta+1)K_1)$, we have
$$\rho( K_{\delta}\nabla (K_{n-(b+1)\delta-1}\cup (b\delta+1)K_1))>\rho(K_{n-b\delta-1})=n-b\delta-2,$$
and hence $\lambda_1(A_\Pi^\delta)>\lambda_1(A_\Pi^s)$. Combining this with Lemma \ref{Q}, (\ref{equ::5}) and (\ref{equ::6}), we have
$$\rho(G)\le\rho(G^1_s)\le\rho(G^2_s)<\rho(K_{\delta}\nabla (K_{n-(b+1)\delta-1}\cup (b\delta+1)K_1)).$$

\noindent{\bf{Case 2.}} $s\le \delta-1$.

Let $G^3_s=K_s\nabla(K_{n-s-(\delta-s+1)(bs+1)}\cup (bs+1)K_{\delta-s+1})$. Recall that $c_i\geq c_{i+1}$ for $1\leq i\leq l-1$.
It is easy to find that $c_l\ge \delta-s+1$ because
the minimum degree of $G^1_s$ is at least $\delta$. By Lemma \ref{Dan Fan}, we have
$$\rho(G^1_s)\le\rho(G^3_s),$$
with equality if and only if $G^1_s\cong G^3_s$.
In what follows, we shall prove that $\rho(G^3_s)<n-b\delta-2$. Note that
$$e(G_3)={n-(\delta+1-s)(bs+1)\choose 2}+(bs+1)\Big(\frac{(\delta+1-s)(\delta-s)}{2}+s(\delta-s+1)\Big).$$
Let
\begin{equation*}
\begin{aligned}
h(n)&=\Big(n-\Big(b+\frac{1}{2}\Big)\delta-\frac{3}{2}\Big)^2-\Big(2e(G_3)-n\delta+\frac{(\delta+1)^2}{4}\Big)\\
    &=2(\delta-s)(bs-b+1)n-(s-1)((\delta+1)s-s^2+\delta)(\delta-s)b^2+(4(\delta+1)s^2-s^3-3(\delta\\
&~~~+1)^2s+a^2+3a)b+(2\delta+2)s-2\delta^2-3\delta.\\
    \end{aligned}
\end{equation*}
As $n\ge b\delta^3+\delta$, $\delta\ge s+1$, $b\ge1$ and $s\ge1$, we have
\begin{equation}\label{equ::9}
\begin{aligned}
h(n)&\ge h(b\delta^3+\delta)\\
&=(s-1)(\delta-s)(2\delta^3-(s+1)\delta+s^2-s)b^2+(2\delta^4-2\delta^3s-(s+1)\delta^2\\
&~~+(2s^2-4s+3)\delta-s^3+4s^2-3s)b+(2\delta+2)s-3\delta+2s.
\end{aligned}
\end{equation}
%If $s=1$, then $$h(n)\ge h(b\delta^3+\delta)=2\delta^2b(\delta^2-\delta-1)+(b-3)+\delta+2>0.$$
Let
\begin{equation*}\label{equ::21}
\begin{aligned}
\omega(\delta)&=2\delta^4-2\delta^3s-(s+1)\delta^2+(2s^2-4s+3)\delta-s^3+4s^2-3s.
\end{aligned}
\end{equation*}
Recall that $\delta\ge s+1$ and $s\geq 1$. Then
\begin{equation*}
\begin{aligned}
\omega(\delta)&\ge 2\delta^3(s+1)-2\delta^3s-(s+1)\delta^2+(2s^2-4s+3)\delta-(\delta-1)^3+4s^2-3s\\
&=\delta^3-(s-2)\delta^2+(2s^2-4s)\delta+4s^2-3s+1\\
&\ge (s+1)\delta^2-(s-2)\delta^2+(2s^2-4s)\delta+4s^2-3s+1\\
&=3\delta^2+2(s^2-2s)\delta+4s^2-3s+1\\
&>\delta.
\end{aligned}
\end{equation*}
Combining this with (\ref{equ::9}), $s\geq 1$, $b\geq 1$ and $\delta\ge s+1$, we can deduce that $h(n)>0$.
It follows that
\begin{equation}\label{equ::11}
\begin{aligned}
n-b\delta-2>\frac{\delta-1}{2}+\sqrt{2e(G_3)-n\delta+\frac{(\delta+1)^2}{4}}.
\end{aligned}
\end{equation}
Since $K_{\delta}\nabla (K_{n-(b+1)\delta-1}\cup (b\delta+1)K_1)$ contains $K_{n-b\delta-1}$ as a proper subgraph, we have
$$\rho(K_{\delta}\nabla (K_{n-(b+1)\delta-1}\cup (b\delta+1)K_1))>\rho(K_{n-b\delta-1})=n-b\delta-2.$$
%We have
%\begin{equation}\label{equ::13}
%\begin{aligned}
%&~~~~n-(b+\frac{1}{2})\delta-\frac{3}{2}-\sqrt{2e(G_3)-n\delta+\frac{(\delta+1)^2}{4}}\\
%         &=n-b\delta-2-\Big(\frac{\delta-1}{2}+\sqrt{2e(G_3)-n\delta+\frac{(\delta+1)^2}{4}}\Big)\\
%         &>0.
%\end{aligned}
%\end{equation}
%That is $$\frac{\delta-1}{2}+\sqrt{2e(G_3)-n\delta+\frac{(\delta+1)^2}{4}}\le n-b\delta-2.$$
According to Lemma \ref{degree} and (\ref{equ::11}), we obtain that
$$\rho(G_3)\leq \frac{\delta-1}{2}+\sqrt{2e(G_3)-n\delta+\frac{(\delta+1)^2}{4}} < n-b\delta-2<\rho(K_{\delta}\nabla (K_{n-(b+1)\delta-1}\cup (b\delta+1)K_1)).$$

\noindent{\bf{Case 3.}} $s=\delta$.

By Lemma \ref{Dan Fan}, we have
$$\rho(G^1_s)\le\rho(K_{\delta}\nabla (K_{n-(b+1)\delta-1}\cup (b\delta+1)K_1)),$$
with equality if and only if $G^1_s\cong K_{\delta}\nabla (K_{n-(b+1)\delta-1}\cup (b\delta+1)K_1)$.
Combining this with (\ref{equ::5}), we deduce that
$$\rho(G)\le\rho(K_{\delta}\nabla (K_{n-(b+1)\delta-1}\cup (b\delta+1)K_1)),$$
with equality if and only if $G\cong K_{\delta}\nabla (K_{n-(b+1)\delta-1}\cup (b\delta+1)K_1)$.

When $\frac{1}{\tau}\geq 1$ is an integer, we can conclude that
$$\rho(G)\le \rho\Big(K_{\delta}\nabla \Big(K_{n-\frac{(\tau+1)\delta}{\tau}-1}\cup \Big(\frac{\delta}{\tau}+1\Big)K_1\Big)\Big),$$
with equality if and only if $G\cong K_{\delta}\nabla \Big(K_{n-\frac{(\tau+1)\delta}{\tau}-1}\cup (\frac{\delta}{\tau}+1)K_1\Big)$. It is easy to verify that $K_{\delta}\nabla \Big(K_{n-\frac{(\tau+1)\delta}{\tau}-1}\cup (\frac{\delta}{\tau}+1)K_1\Big)$ is not $\tau$-tough. This completes the proof of (ii).\end{proof}

\section{Proof of Theorem \ref{td}}
\begin{lemma}[See \cite{Zhai-Lin}]\label{Sum}
Suppose that $G$ is a connected graph and $X=\left(x_1, x_2, \ldots, x_n\right)^{T}$ is the Perron vector of $G$. Let $S_1, S_2, T$ be three mutually disjoint and non-empty subsets of $V(G)$ with $e\left(T, S_1\right)=0$ and $e(T, S_2)=|T||S_2|$, and let $G'=G+\{i j \mid i \in S_1, j \in T\}-\{ij \mid i \in S_2, j \in T\}$. If $\sum_{i \in S_1} x_i \geq \sum_{i \in S_2} x_i$, then $\rho(G')>\rho(G)$.
\end{lemma}
\begin{lemma}[See \cite{Nosal}]\label{H}
Let $G$ be a bipartite graph with $m$ edges. Then
$$\rho(G)\le\sqrt{m},$$
with equality if and only if $G\cong K_{p,q}\cup(n-p-q)K_1$, where $pq=m$.
\end{lemma}

\begin{lemma}\label{H2B2}
Let $k\geq 2$ be an integer and $n\geq 2k^2+6k$ be an even number. If $2k\mid n$, then
$$\rho(K_{\frac{n}{2}-1,\frac{n}{2}-\frac{n}{2k}}\nabla_1 O_{1,\frac{n}{2k}})> \rho(K_{k-1,\frac{n}{2}-1}\nabla_1 O_{\frac{n}{2}-k+1,1}).$$
\end{lemma}

\renewcommand\proofname{\bf Proof}
\begin{proof}
Let $G=K_{k-1,\frac{n}{2}-1}\nabla_1 O_{\frac{n}{2}-k+1,1}$. Then $A(G)$ has the equitable quotient matrix
\begin{align*}
B_{\Pi_1}=\left[
\begin{matrix}
  0 & \frac{n}{2}-1 & 0 & 1 \\
  k-1 & 0 & \frac{n}{2}-k+1 & 0 \\
  0 & \frac{n}{2}-1 & 0 & 0  \\
  k-1 & 0 & 0 & 0
\end{matrix}
  \right].
\end{align*}
By a simple calculation, the characteristic polynomial of $B_{\Pi_1}$ is
\begin{equation*}
\mu_1(B_{\Pi_1},x)=x^4-\Big(k-1+\frac{1}{4}n^2-\frac{1}{2}n\Big)x^2-\frac{1}{4}(n-2)(k-1)(2k-n-2).
\end{equation*}
Take $G^*=K_{\frac{n}{2}-1,\frac{n}{2}-\frac{n}{2k}}\nabla_1 O_{1,\frac{n}{2k}}$. Thus, $A(G^*)$ has the equitable quotient matrix
\begin{align*}
C_{\Pi_2}=\left[
\begin{matrix}
  0 & \frac{n}{2}-\frac{n}{2k} & 0 & \frac{n}{2k} \\
  \frac{n}{2}-1 & 0 & 1 & 0 \\
  0 & \frac{n}{2}-\frac{n}{2k} & 0 & 0  \\
  \frac{n}{2}-1 & 0 & 0 & 0
\end{matrix}
  \right].
\end{align*}
The characteristic polynomial of $C_{\Pi_2}$ is
$$\mu_2(C_{\Pi_2},x)=x^4-\frac{1}{4k}n(kn-2)x^2+\frac{1}{8k^2}n^2(k-1)(n-2).$$
For $x> \sqrt{\frac{n}{2}(\frac{n}{2}-1)}$, we have
\begin{equation*}
\begin{aligned}
&\mu_1(B_{\Pi_1},x)-\mu_2(C_{\Pi_2},x)\\
&=\frac{(k-1)(n-2k)x^2}{2k}+\frac{(n-2)(n-2k)(k-1)(2k^2-2k-n)}{8k^2}\\
&> \frac{1}{8k^2}(k-1)^2(n-2k)(n+2k)(n-2)
~~~(\mbox{since $x>\sqrt{ \frac{n}{2}(\frac{n}{2}-1)}$})\\
&> 0~~~(\mbox{since $k\ge2$, $n\ge2k^2+6k$}).
\end{aligned}
\end{equation*}
Therefore, $\mu_1(B_{\Pi_1},x)>\mu_2(C_{\Pi_2},x)$ for $x> \sqrt{\frac{n}{2}(\frac{n}{2}-1)}$.
Since $K_{\frac{n}{2}-1,\frac{n}{2}-\frac{n}{2k}}\nabla_1 O_{1,\frac{n}{2k}}$ contains $K_{\frac{n}{2}-1,\frac{n}{2}}$ as a proper subgraph,
$\rho(K_{\frac{n}{2}-1,\frac{n}{2}-\frac{n}{2k}}\nabla_1 O_{1,\frac{n}{2k}})>\sqrt{\frac{n}{2}(\frac{n}{2}-1)}$ by Lemma \ref{edge}, and hence  $\lambda_1(C_{\Pi_2})>\lambda_1(B_{\Pi_1})$. According to Lemma \ref{Q}, we have
$$\rho(K_{\frac{n}{2}-1,\frac{n}{2}-\frac{n}{2k}}\nabla_1 O_{1,\frac{n}{2k}})>\rho(K_{k-1,\frac{n}{2}-1}\nabla_1 O_{\frac{n}{2}-k+1,1}),$$
as required.
\end{proof}

\begin{lemma}\label{S}
If $1\le s\le \frac{n-4}{4}$, then
$\rho(K_{s,\frac{n}{2}-s-1}\nabla_{1}O_{\frac{n}{2}-s,s+1})>\rho(K_{s+1,\frac{n}{2}-s-2}\nabla_{1}O_{\frac{n}{2}-s-1,s+2}).$
\end{lemma}
\renewcommand\proofname{\bf Proof}
\begin{proof}
Let $G^s=K_{s,\frac{n}{2}-s-1}\nabla_{1}O_{\frac{n}{2}-s,s+1}$. Then $A(G^s)$ has the equitable quotient matrix
\begin{align*}
M^s_{\Pi_4}=\left[
\begin{matrix}
  0 & \frac{n}{2}-s-1 & 0 & s+1 \\
  s & 0 & \frac{n}{2}-s & 0 \\
  0 & \frac{n}{2}-s-1 & 0 & 0  \\
  s & 0 & 0 & 0
\end{matrix}
  \right].
\end{align*}
By a simple calculation, the characteristic polynomial of $M^s_{\Pi_4}$ is
$$\phi(M^s_{\Pi_4},x)=x^4-\Big(s^2+s+\frac{1}{4}n^2-\frac{1}{2}ns-\frac{1}{2}n\Big)x^2+\frac{1}{4}s(s+1)(n-2s)(n-2s-2).$$
Assume that $G^{s+1}=K_{s+1,\frac{n}{2}-s-2}\nabla_{1}O_{\frac{n}{2}-s-1,s+2}$. Thus, $A(G^{s+1})$ has the equitable quotient matrix $M^{s+1}_{\Pi_4}$, which is obtained by
replacing $s$ with $s+1$ in $M^s_{\Pi_4}$. For $x>\sqrt{ \frac{sn}{2}}$, we have
\begin{equation*}
\begin{aligned}
&\phi(M^{s+1}_{\Pi_4},x)-\phi(M^{s}_{\Pi_4},x)\\
&=\frac{1}{2}(n-4s-4)(x^2+n(s+1)-2s^2-4s-2)\\
&>\frac{1}{2}(n-4s-4)\Big(\Big(\frac{3s}{2}+1\Big)n-2s^2-4s-2\Big)~~~(\mbox{since $x>\sqrt{ \frac{sn}{2}}$})\\
&\geq 0 ~(\mbox{since $n\geq 4s+4$ and $s\geq 1$}).
\end{aligned}
\end{equation*}
Thus, $\phi(M^{s+1}_{\Pi_4},x)>\phi(M^{s}_{\Pi_4},x)$ for $x>\sqrt{ \frac{sn}{2}}$.  Since $K_{s,\frac{n}{2}}$ is a proper subgraph of $G^{s}$,  $\rho(G^{s})>\rho(K_{s,\frac{n}{2}})=\sqrt{\frac{sn}{2}}$.  This implies that $\lambda_1(M^{s}_{\Pi_4})>\lambda_1(M^{s+1}_{\Pi_4})$. By Lemma \ref{Q}, we get $\rho(G^s)>\rho(G^{s+1})$, as required.\end{proof}

Let $\mathcal{G}_{n,r}$ be a set of connected balanced bipartite graph of even order $n$ which are not $r$-tough, where $r\geq 2$ or $\frac{1}{r}$ is a positive integer.
\renewcommand\proofname{\bf Proof of Theorem \ref{td}}
\begin{proof}
Suppose that $G'\in \mathcal{G}_{n,r}$ attains the maximum spectral radius.
Thus, for $G\in \mathcal{G}_{n,r}$, we can obtain that
\begin{equation}\label{equ::12}
\begin{aligned}
\rho(G)\leq \rho(G').
\end{aligned}
\end{equation}
Since $G'=(X,Y)$ is not $r$-tough, without loss of generality, we assume that there exists some nonempty proper subset $S\subset X$ such that $|S|<r(c(G'-S)-1)$ with $c(G'-S)>1$. Take $|S|=s$ and $c(G'-S)=q$. Then $s\le r(q-1)-1$.
Suppose that $G_1, G_2,\cdots, G_t$ are the non-trivial components of $G'-S$.
For $1\leq i\leq t$, let $G_i=(X_i,Y_i)$ with $X_i\subseteq X$ and $Y_i\subseteq Y$, and let $|X_i|=n_i$ and $|Y_i|=m_i$. By the maximality of $\rho(G')$, we can deduce that $G_i\cong K_{n_i,m_i}$ for $1\leq i\leq t$ and $G[S\cup Y]\cong K_{s,\frac{n}{2}}$. We first claim that $t=1$. If not, $t\geq 2$. Without loss of generality, we assume that $\sum_{v\in Y_1}x(v)=\max\{\sum_{v\in Y_i}x(v)|~1\leq i\leq t\}$.
Let $G''=G'-\{uv|~u\in X_2, v\in Y_2\}+\{uv|~u\in X_2, v\in Y_1\}$.
By Lemma \ref{Sum}, we have $\rho(G'')>\rho(G')$, which contradicts the maximality of $\rho(G')$. Hence, $t=1$. It follows that $G'\cong K_{s,\frac{n}{2}-q+1}\nabla_1 O_{\frac{n}{2}-s,q-1}$.

(i) Let $r\geq 2$ be an integer. Since $r(q-1)-1\le \frac{n}{2}-1$, we have $2\le q\le\frac{n}{2r}+1$. Note that
\begin{equation*}
\begin{aligned}
e(G')&=\frac{n}{2}(r(q-1)-1)+\Big(\frac{n}{2}-r(q-1)+1\Big)\Big(\frac{n}{2}-q+1\Big)\\
&=rq^2-\Big(\frac{n}{2}+2r+1\Big)q+\frac{n^2}{4}+\frac{n}{2}+r+1\\
&=\omega(q).
\end{aligned}
\end{equation*}
If $3\le q\le \frac{n}{2r}$, then
$$\omega\Big(\frac{n}{2r}\Big)-\omega(3)=\Big(\frac{1}{2}-\frac{1}{2r}\Big)n-3r+3\geq r^2-r> 0$$
due to $n\ge2r^2+6r$ and $r\geq 2$. This implies that the maximum value of $\omega(q)$ is attained at $q=\frac{n}{2r}$ for $3\le q\le \frac{n}{2r}$.
Since
$$\omega\Big(\frac{n}{2r}\Big)=\frac{1}{4}n^2-\Big(\frac{1}{2}+\frac{1}{2r}\Big)n+r+1=\frac{n}{2}\Big(\frac{n}{2}-1\Big)-\Big(\frac{n}{2r}-r-1\Big)<\frac{n}{2}\Big(\frac{n}{2}-1\Big),$$ due to
$n\ge2r^2+6r$, we have $\omega(q)\leq \omega(\frac{n}{2r})<\frac{n}{2}(\frac{n}{2}-1)$. Combining this with Lemma \ref{H}, we obtain that
\begin{equation}\label{equ::13}
\begin{aligned}
\rho(G')\leq \sqrt{e(G')}<\sqrt{\frac{n}{2}\Big(\frac{n}{2}-1\Big)}.
\end{aligned}
\end{equation}
Note that $K_{r-1,\frac{n}{2}-1}\nabla_1 O_{\frac{n}{2}-r+1,1}\in \mathcal{G}_{n,r}$ contains $K_{\frac{n}{2},\frac{n}{2}-1}$ as a proper subgraph. Then
$$\sqrt{\frac{n}{2}\Big(\frac{n}{2}-1\Big)}<\rho(K_{r-1,\frac{n}{2}-1}\nabla_1 O_{\frac{n}{2}-r+1,1})\le\rho(G'),$$
which contradicts (\ref{equ::13}). If $q=2$, then $s\leq r(q-1)-1=r-1$, and hence
$G'\cong K_{r-1,\frac{n}{2}-1}\nabla_1 O_{\frac{n}{2}-r+1,1}.$ Combining this with (\ref{equ::12}), we have
\begin{equation}\label{equ::14}
\begin{aligned}
\rho(G)\leq \rho(K_{r-1,\frac{n}{2}-1}\nabla_1 O_{\frac{n}{2}-r+1,1}),
\end{aligned}
\end{equation}
with equality if and only if $G\cong K_{r-1,\frac{n}{2}-1}\nabla_1 O_{\frac{n}{2}-r+1,1}$. Next, we consider the case of $q=\frac{n}{2r}+1$. If $2r\mid n $, then $s\leq r(q-1)-1=\frac{n}{2}-1$. By the maximality of $\rho(G')$, we have $G'\cong K_{\frac{n}{2}-1,\frac{n}{2}-\frac{n}{2r}}\nabla_1 O_{1,\frac{n}{2r}}$. Combining this with (\ref{equ::12}), (\ref{equ::14}) and Lemma \ref{H2B2}, we deduce that
$$\rho(G)\leq \rho(K_{\frac{n}{2}-1,\frac{n}{2}-\frac{n}{2r}}\nabla_1 O_{1,\frac{n}{2r}}),$$
with equality if and only if $G\cong K_{\frac{n}{2}-1,\frac{n}{2}-\frac{n}{2r}}\nabla_1 O_{1,\frac{n}{2r}}$. Notice that $K_{\frac{n}{2}-1,\frac{n}{2}-\frac{n}{2r}}\nabla_1 O_{1,\frac{n}{2r}}\in \mathcal{G}_{n,r}$. Thus, the result as follows. If $2r \nmid n$, then $q=\lfloor\frac{n}{2r}\rfloor+1$, and hence
$$s\leq r(q-1)-1=r\Big\lfloor\frac{n}{2r}\Big\rfloor-1.$$
Again by the maximality of $\rho(G')$, we get $G'\cong K_{r\lfloor\frac{n}{2r}\rfloor-1,\frac{n}{2}-\lfloor\frac{n}{2r}\rfloor}\nabla_1 O_{\frac{n}{2}-r\lfloor\frac{n}{2r}\rfloor+1,\lfloor\frac{n}{2r}\rfloor}$.
 Combining this with (\ref{equ::12}) and (\ref{equ::14}), we deduce that
$$\rho(G)\leq \max\{\rho(K_{r-1,\frac{n}{2}-1}\nabla_1 O_{\frac{n}{2}-r+1,1}),\rho(K_{r\lfloor\frac{n}{2r}\rfloor-1,\frac{n}{2}-\lfloor\frac{n}{2r}\rfloor}\nabla_1 O_{\frac{n}{2}-r\lfloor\frac{n}{2r}\rfloor+1,\lfloor\frac{n}{2r}\rfloor})\},$$
where the equality holds if and only if $G\cong K_{r-1,\frac{n}{2}-1}\nabla_1 O_{\frac{n}{2}-r+1,1}$ or $G\cong  K_{r\lfloor\frac{n}{2r}\rfloor-1,\frac{n}{2}-\lfloor\frac{n}{2r}\rfloor}\nabla_1 $ $O_{\frac{n}{2}-r\lfloor\frac{n}{2r}\rfloor+1,\lfloor\frac{n}{2r}\rfloor}$.
It is easy to find that $K_{r\lfloor\frac{n}{2r}\rfloor-1,\frac{n}{2}-\lfloor\frac{n}{2r}\rfloor}\nabla_1$ $O_{\frac{n}{2}-r\lfloor\frac{n}{2r}\rfloor+1,\lfloor\frac{n}{2r}\rfloor}\in \mathcal{G}_{n,r}$. This completes the proof of (i).

(ii)
Let $b=\frac{1}{r}$ be a positive integer. Then $q\geq bs+2$. By the maximality of $\rho(G')$, we get $G'\cong K_{s,\frac{n}{2}-bs-1}\nabla_{1}O_{\frac{n}{2}-s,bs+1}$.
Now, we consider the following two case depending on the value of $b$.

{\flushleft\bf Case 1.} $b=1$.

In this case, we may assume that $s\leq \frac{n}{2}-s-1$. Thus, $1\leq s\leq \frac{n-2}{4}$. According to Lemma \ref{S}, we conclude that $G'\cong K_{1,\frac{n}{2}-2}\nabla_{1}O_{\frac{n}{2}-1,2}$. Combining this with (\ref{equ::12}), we have
$$\rho(G)\leq \rho(K_{1,\frac{n}{2}-2}\nabla_{1}O_{\frac{n}{2}-1,2}),$$
with equality if and only if $G\cong K_{1,\frac{n}{2}-2}\nabla_{1}O_{\frac{n}{2}-1,2}$.

{\flushleft\bf Case 2.} $b\geq 2$.

As $G'$ is connected, we have $\frac{n}{2}-bs-1\ge1$. Thus, $1\le s\le \frac{n-4}{2b}$.
If $s=1$, then $G'\cong K_{1,\frac{n}{2}-b-1}\nabla_{1}O_{\frac{n}{2}-1,b+1}$. Furthermore, by (\ref{equ::12}), we obtain that
\begin{equation*}
\begin{aligned}
\rho(G)\leq\rho(K_{1,\frac{n}{2}-b-1}\nabla_{1}O_{\frac{n}{2}-1,b+1}),
\end{aligned}
\end{equation*}
with equality if and only if $G\cong K_{1,\frac{n}{2}-b-1}\nabla_1 O_{\frac{n}{2}-1,b+1}$.
For $2\leq s\leq \frac{n-4}{2b}$, we have
\begin{equation}\label{equ::15}
\begin{aligned}
\rho(G')&\le\sqrt{e(G')}\\
&=\sqrt{\frac{ns}{2}+\Big(\frac{n}{2}-s\Big)\Big(\frac{n}{2}-bs-1\Big)}\\
&=\sqrt{bs^2-\Big(\frac{bn}{2}-1\Big)s+\frac{n^2}{4}-\frac{n}{2}}
\end{aligned}
\end{equation}
by Lemma \ref{H}. Let $\gamma(s)=bs^2-(\frac{bn}{2}-1)s+\frac{n^2}{4}-\frac{n}{2}$. By a simple calculation, we have
$$\gamma(2)-\gamma\Big(\frac{n-4}{2b}\Big)=\frac{1}{4b}((b-1)n-4b+2)(n-4b-4)> 0,$$
where the inequality follows from the fact that $n\ge 4b+4$ and $b\ge2$. This implies that
for $2\le s\le \frac{n-4}{2b}$, the maximum value of $\gamma(s)$ is attained at $s=2$. Combining this with
(\ref{equ::15}), we deduce that
\begin{equation}\label{equ::16}
\begin{aligned}
\rho(G')&\le\sqrt{\gamma(2)}\\
&=\sqrt{\frac{1}{4}n^2-\Big(b+\frac{1}{2}\Big)n+4b+2}\\
&=\sqrt{\frac{n}{2}\Big(\frac{n}{2}-1-b\Big)-\Big(\frac{1}{2}n-4\Big)b+2}\\
&<\sqrt{\frac{n}{2}\Big(\frac{n}{2}-1-b\Big)} ~(\mbox{since $n\geq 4b+4$ and $b\ge2$)}.
\end{aligned}
\end{equation}
Since
 $K_{1,\frac{n}{2}-b-1}\nabla_1 O_{\frac{n}{2}-1,b+1}$ contains $K_{\frac{n}{2},\frac{n}{2}-b-1}$ as a proper subgraph, we have
\begin{equation*}
\begin{aligned} \sqrt{\frac{n}{2}\Big(\frac{n}{2}-1-b\Big)}=\rho(K_{\frac{n}{2},\frac{n}{2}-b-1})<\rho(K_{1,\frac{n}{2}-b-1}\nabla_1 O_{\frac{n}{2}-1,b+1})\le \rho(G'),
\end{aligned}
\end{equation*}
which contradicts  (\ref{equ::16}).

For the positive integer $\frac{1}{r}$, by the discussions as above, we have
$$\rho(G)\leq\rho(K_{1,\frac{n}{2}-\frac{1}{r}-1}\nabla_{1}O_{\frac{n}{2}-1,\frac{1}{r}+1}),$$
with equality if and only if $G\cong K_{1,\frac{n}{2}-\frac{1}{r}-1}\nabla_1 O_{\frac{n}{2}-1,\frac{1}{r}+1}$. Moreover, it is easy to verify that $K_{1,\frac{n}{2}-\frac{1}{r}-1}\nabla_1 O_{\frac{n}{2}-1,\frac{1}{r}+1}\in \mathcal{G}_{n,r}$ where $\frac{1}{r}$ is a positive integer. Thus, the result follows.\end{proof}

\vskip 0.382 in
\noindent
\textbf{\Large Declaration of Competing Interest}

The authors declare that they have no known competing financial interests or personal
relationships that could have appeared to influence the work reported in this paper.

%\begin{prop}
%Let $G$ and $\overline{G}$ are both connected and $1/2<\alpha<1$. If $G$ have pendent vertices, then
%$$\lambda_2(A_\alpha(G))+\lambda_2(A_\alpha(\overline{G}))\le(2n-3)\alpha-1.$$
%\end{prop}
%\begin{proof}
%Since $G$ and $\overline{G}$ are both connected,
%$\Delta(G)<n-1$ and $\Delta(\overline{G})<n-1$. We get that
%$\lambda_2(A_\alpha(G))\le (n-2)\alpha$ and $\lambda_2(A_\alpha(\overline{G}))<(n-2)\alpha$
%by Lemma \ref{k+1}.
%Note that $\lambda_2(A_\alpha(G^*))=(n-1)\alpha-1$, we have
%$\lambda_2(A_\alpha(G))\le\lambda_2(A_\alpha(G^*))=(n-1)\alpha-1$ by Lemma \ref{edge}.
%Since the equalities can not both holds,
%$$\lambda_2(A_\alpha(G))+\lambda_2(A_\alpha(\overline{G})<(n-1)\alpha-1+(n-2)\alpha=(2n-3)\alpha-1.$$
%\end{proof}

\end{document}